\documentclass[12pt]{article}

\begin{document}

\author{Igor Borovikov}
\title{Order Reduction of Optimal Control Systems}
\date{September 6, 2004}
\maketitle

\begin{abstract}
The paper presents necessary and sufficient conditions for the order
reduction of optimal control systems. Exploring the corresponding
Hamiltonian system allows to solve the order reduction problem in terms of
dynamical systems, observability and invariant differential forms. The
approach is applicable to non-degenerate optimal control systems with smooth
integral cost function. The cost function is defined on the trajectories of
a smooth dynamical control system with unconstrained controls and fixed
boundary conditions. Such systems form a category of Lagrangian systems with
morphisms defined as mappings preserving extremality of the trajectories.
Order reduction is defined as a factorization in the category of Lagrangian
systems.
\end{abstract}



\section{Introduction}

The invariance in nonlinear control theory allows to approach and solve many
important problems in an effective way. In particular, the invariance of
dynamical control system with respect to the action of Lie group or Lie
algebra allows order reduction by factoring out a system of smaller
dimension. The review \cite{pa1} and monographs \cite{nvs}, \cite{el}
present the current state of the art in this area.

The study of invariance in optimal control (like in classical works \cite{vs}%
, \cite{su} or in recent publications \cite{ech}, \cite{to1}) can be built
on the results of classical mechanics \cite{ab}, \cite{ar}. This places the
geometric approach to optimal control into the rich context of Hamiltonian
mechanics and symplectic geometry.

Different generalizations of the symplectic geometry allow to study wider
range of symmetries of optimal control systems. In the recent work \cite{ech}
an optimal control system is treated as a Hamiltonian system on the
corresponding presymplectic manifold. The symmetries of optimal control
system are symplectic (when considered on a special symplectic subspace)
actions of Lie group that leave both dynamical system and Lagrange function
invariant. It is shown that such symmetries allow to reduce the order of the
optimal control system.

A more general point of view, based on theory of categories and theory of
decomposition, allows to formulate the problem of order reduction in terms
of general factor- and sub- objects (see \cite{pa1} for review of theory of
decomposition for control systems).

In \cite{ch} the order reduction of a smooth variational problem

$$
\int_{0}^{T}L(q,\dot q)dt \rightarrow extr, \qquad q \in R^n, q(0)=q_0,
q(T)=q_1 
$$

was studied in categorial framework. Variational systems form a category
with morphisms preserving extremals. The extremals are solutions of the
corresponding Euler-Lagrange equations. These equations form a dynamical
system, which can be studied and reduced using general geometric methods.
Solving the inverse problem of variational calculus allows to rebuild the
reduced variational system, when it exists, from the corresponding dynamical
factorsystem. The work \cite{ch} presents necessary and sufficient
conditions for order reduction (factorization) of variational systems.

One of the advantages of the categorial framework is that it allows to pose
and solve the order reduction problem in the most general and complete way:
the results of \cite{ch} cover {\it all} factorizations in the category of
variational systems with morphisms preserving extremals.

Here we will define a more general category of Lagrangian systems that
includes variational systems from \cite{ch} as a sub-category. The objects
of the new category are optimal control systems with smooth integral
criterion, smooth dynamical control system and unconstrained controls.
Morphisms in the category of Lagrangian systems are also defined as mappings
preserving extremals. Pontryagin's maximum principle allows to convert
Lagrangian system into corresponding Hamiltonian system. We will define a
new category of Hamiltonian systems with morphisms that may not necessarily
preserve underlying symplectic structure. That will allow to study
factorizations of Hamiltonian systems using methods for general dynamical
systems. Finally we will show that factorizations of Hamiltonian systems and
Lagrangian systems are corresponding to each other provided that
observability condition for Hamiltonian factorsystem is satisfied.

The main results of this work were announced in \cite{bo1} and \cite{bo2}.


\section{Order Reduction as Factorization}

The problem of order reduction could be approached for each particular type
of mathematical objects individually. However we can study it in a more
uniform way by placing it into the categorial framework where order
reduction can be treated as factorization in the appropriate category.

For a brief informal illustration we will consider factorization of smooth
dynamical systems. Smooth dynamical systems will form a category $DS$ if we
define morphisms - mappings of the dynamical systems into each other. A
natural morphism is a diffeomorphisms from one dynamical system into the
other dynamical system that maps trajectories of the original dynamical
system into the trajectories of the image system. Factorization is a special
morphism that does not introduce anything additional to the factor object,
which could not be derived from the original object. A morphism in category $%
DS$ defines a factorization if it is a surjective submersion (mapping onto
of the full rank). A factor system in this case is the dynamical system on
the image (factor) space. Any trajectory of the factor system has its
inverse image - a set of trajectories of the original system. Also all
trajectories of the original system map into some trajectories of the factor
system. It is well known that such systems are described by $f$-related
vector fields and provide the classical example of order reduction. Thus the
order reduction of dynamical systems can be naturally described in very
general terms of factorization.

A formal definitions for this approach could be derived from theory of
structures \cite{bo} or, equivalently, from theory of categories \cite{bu}.
For in depth discussion of theory of decomposition in application to control
systems see \cite{pa1}.

In this paper we will apply categorial approach to the order reduction
problem for optimal control systems.

\section{The Category of Lagrangian Systems}

We will consider optimal control system formed by a smooth control system
and an integral cost function:

\begin{eqnarray}
\int_{0}^{T}L(q,u)dt \rightarrow extr  \label{L1} \\
\dot q = f(q,u)  \label{UDS1}
\end{eqnarray}

where $q$ is an $n$-dimensional vector of phase variables and $u$ is an $m$%
-dimensional vector of controls. The vector of control $u$ is not
constrained, so $u\in R^m$. We will call such systems {\it Lagrangian systems%
}, because (under the appropriate conditions) the optimal control problem (%
\ref{L1},\ref{UDS1}) with fixed boundary conditions is equivalent to
Lagrange variational problem. The equivalent Lagrangian variational problem
can be obtained by eliminating controls $u$ from both (\ref{L1}) and (\ref
{UDS1}), then transforming \ref{UDS1} to the implicit form $F(q,\dot q)=0$.
Detailed definitions follow later.

We will assume that the optimal control system (\ref{L1}, \ref{UDS1})
defines a field of extremal solutions. An extremal solution for (\ref{L1}, 
\ref{UDS1}) is a curve $\gamma (t) = (q(t), u(t))$ such that it satisfies (%
\ref{UDS1}) and the functional (\ref{L1}) achieves on $\gamma(t)$ an
extremal value within the class of the curves with fixed boundary points. As
usually, extremals are not necessarily optimal curves: extremality means
that the curve has vanishing conditional variations of the functional in its
vicinity.

We will be looking for an optimal control system of the same kind but of a
lower order (a problem of order reduction or {\it factorization}):

\begin{eqnarray}
\int_{0}^{T} Q(y,v)dt \rightarrow extr  \label{L2} \\
\dot y = F(y,v)  \label{UDS2}
\end{eqnarray}

such that there exists mapping $y= y(q,u), v=v(q,u)$, which maps extremals
of (\ref{L1}, \ref{UDS1}) to the extremals of (\ref{L2}, \ref{UDS2}).

In this paper we will derive the necessary and sufficient conditions for the
factorization of optimal control systems (\ref{L1}, \ref{UDS1}). We will use
Pontryagin's Maximum Principle to transform the optimal systems into
Hamiltonian systems. Hamiltonian systems will allow for an intuitive
geometric approach to the problem of factorization. Finally we will
translate the results obtained in the terms of Hamiltonian systems back into
the domain of the original optimal control systems.

Everywhere we assume smoothness an locality: all manifolds are open simply
connected regions of $R^n$, all functions are smooth (have as many
derivatives as necessary), so we drop adjectives ''smooth'' and ''local'' in
most cases. We use Einstien's convention for summation: terms with repeating
subscript and superscript index automatically sums.

\section{Factorization of Hamiltonian systems}

A Hamiltonian system is a dynamical system generated by the gradient flow of
a Hamiltonian function defined on a symplectic manifold. More formally, a
triplet $HS=(M,\omega ^2,H)$ defines a Hamiltonian system on a symplectic
manifold $M,dimM=2n$, with the symplectic structure $\omega ^2$ and
Hamiltonian function $H:M\to R$.

A Hamiltonian system $HS$ defines canonical equations in the form:

\begin{eqnarray}
\dot z = IdH(z), \qquad z \in M  \label{HS0}
\end{eqnarray}

where $I:TM\to T^{*}M$ is the isomorphism induced by the symplectic
structure. We will also use $HS$ to denote the canonical equations (\ref{HS0}%
).

In local canonical coordinates $(p,q)$ on $M$, we have:

\begin{eqnarray}
\omega ^2 = dp_i \wedge dq^i, \qquad i = \overline{1,n}  \label{omega2}
\end{eqnarray}

and the system (\ref{HS0}) has the form

$$
\dot p = - \partial H / \partial q, \qquad \dot q = \partial H / \partial p 
$$

Let $HS^{\prime}=(N, \overline{\omega}^2, G), \: dim N = 2n$ be another
Hamiltonian system and let $\phi: M \to N$ be a smooth (not necessarily
symplectic) mapping. If $(p,q)$ and $(x,y)$ are canonical coordinates on $M$
and $N$, then we can write $\phi$ in coordinates as $x=x(p,q), y=y(p,q)$.

\newtheorem{definition}{Definition}

\begin{definition}
\label{definition:morphism1} Mapping $\phi: M \to N$ is called a morphism of
Hamiltonian systems if for any solution $z(t), t \in [0,T]$ of the system $HS
$ its image under the mapping $\phi$ is a solution of the system $HS^{\prime}
$ on the interval $t \in [0,T]$.
\end{definition}

Using categorial notation we will write $\phi :HS\to HS^{\prime }$ for the
morphism $\phi $ of Hamiltonian systems $HS$ and $HS^{\prime }$.

The underlying symplectic structure does not participate in the definition
of the morphism directly. So, if $\phi :HS\to HS^{\prime }$ is a morphism,
then by ''erasing'' symplectic structures $\omega ^2$ and $\overline{\omega }%
^2$ from $HS$ and $HS^{\prime }$ respectively we will obtain a morphism of
general dynamical systems. This observation leads to:

\newtheorem{proposition}{Proposition}

\begin{proposition}
\label{proposition:fieldsconnection} Let $HS=(M,\omega ^2, H)$ and $%
HS^{\prime}=(N, \overline{\omega}^2, G)$ be Hamiltonian systems, and let $%
\phi: M \to N$ be a smooth mapping, then the following conditions are
equivalent:

\begin{enumerate}
\item  $\phi$ is morphism of Hamiltonian systems $HS$ and $HS^{\prime}$

\item  vector fields of the systems $HS$ and $HS^{\prime}$ are $\phi$-related
\end{enumerate}
\end{proposition}

If we denote by $I^{\prime}$ the natural isomorphism $T^*N \to TN$ induced
by $\overline{\omega}^2$, then the relation between vector fields of the
Hamiltonian systems $HS$ and $HS^{\prime}$ will have the form

$$
I^{\prime}dG = \phi _* IdH 
$$

Using Poisson bracket $()_{p,q}$ on $M$ we can write the same relation in
canonical coordinates on $M$ as

\begin{eqnarray}
\partial G / \partial x_i = (H, y^i)_{p,q}, \qquad - \partial G / \partial
y^i = (H, x_i)_{p,q}  \label{dG}
\end{eqnarray}

We will focus on the case when $dimN<dimM$ and the mapping $\phi $ is onto
and of full rank, which corresponds to a factorization of Hamiltonian
systems. Let's find when $\phi :M\to N$ maps a Hamiltonian system $HS$ from $%
(M,\omega ^2)$ into a Hamiltonian system on $(N,\overline{\omega }^2)$. The
symplectic form $\overline{\omega }^2$ on the factor space $(N,\overline{%
\omega }^2)$ has to be an invariant of the factor-system by the definition
of the Hamiltonian system. The form $\Omega ^2=\phi ^{*}\overline{\omega }^2$
on $M$ is induced by $\phi $ from the form $\overline{\omega }^2$. If the
canonical coordinates $(p,q)$ and $(x,y)$ on both spaces are fixed, then $%
\phi $ has the form $x=x(p,q),y=y(p,q)$, and we obtain a coordinate
representation

$$
\Omega ^ 2 = dx_i(p,q) \wedge dy^i(p,q), \qquad i=\overline{1,m}. 
$$

\begin{proposition}
\label{proposition:HS_F_system} Let $\phi: M \to N$ be a surjective
submersion. The projection $\phi_*IdH$ of the field $IdH$ from $M$ to $N$
exists and is a Hamiltonian vector field on $(N, \overline{\omega}^2)$ iff
1-form $i_{IdH}\Omega^2$ is closed.
\end{proposition}

Here $i_ab$ denotes the internal product of a vector field $a$ and a form $b$%
.

Due to the locality, a closed form is automatically exact. That means it
actually is a differential of some function on $M$ and, as we will show, on $%
N$ as well.

Note again that the fact that the vector field on $M$ is Hamiltonian is not
used anywhere.

Proof. Necessity. By assumption $\phi_*IdH = I^{\prime}dG$, where $%
I^{\prime}dG$ is a Hamiltonian field on $(N, \overline{\omega}^2)$. Consider
function $\overline{G}=G \circ \phi$, which is defined on $M$. The following
is valid for $\overline{G}$:

$$
d \overline{ G } = d(G \circ \phi) = \phi^* dG = \phi^* i_{I^{\prime}dG} 
\overline{\omega}^2 = i_{IdH} \phi^*\overline{\omega}^2 = i_{IdH} \Omega^2, 
$$

These equalities follow from the chain rule applied to $G\circ \phi $ and
the equivalent transformation of a gradient 1-form into a vector field on a
symplectic manifold.

To receive the same result in coordinates, let's unfold $i_{IdH}\Omega ^2$
in the coordinates $(p,q)$ on $M$:

$$
i_{IdH}\Omega ^2\equiv i_{IdH}dx_i(p,q)\wedge dy^i(p,q) 
$$

Calculating the inner product in the right hand part we will get:

$$
i_{IdH} dx_i(p,q) \wedge dy^i(p,q) = (y^i(p,q),H)_{p,q}dx_i(p,q) -
(x_i(p,q),H)_{p,q}dy^i(p,q) 
$$

The right hand part here is the full differential of $G(x(p,q),y(p,q))$
because of $\phi$-relation of the vector fields expressed by (\ref{dG}).

Sufficiency.

The equality $d\overline{G}=i_{IdH}\Omega ^2$ is equivalent to the fact that
the gradient $d\overline{G}$ can be linearly combined from the gradients of
the independent mapping functions $x(p,q)$ and $y(p,q)$. In coordinates:

\begin{eqnarray}
\frac{\partial{\overline{G}}}{\partial{p_k}} = (y^i,H)_{(p,q)} \frac{\partial{x_i}}{\partial{p_k}} - (x_j,H)_{(p,q)} \frac{\partial{y^j}}{\partial{p_k}}
\qquad  \label{dGdp} \\
\frac{\partial{\overline{G}}}{\partial{q^k}} = (y^i,H)_{(p,q)} \frac{\partial{x_i}}{\partial{q^k}} - (x_j,H)_{(p,q)} \frac{\partial{y^j}}{\partial{q^k}}
\qquad  \label{dGdq}
\end{eqnarray}

\begin{flushleft}
where $k=\overline{1,n}$.
\end{flushleft}

Because of the linear dependence of the gradients we have $\overline{G}%
=G\circ \phi $ with some function $G:N\to R$. Hence $d\overline{G}=\phi
^{*}dG$ and we have a new equality $\phi ^{*}dG=i_{IdH}\phi ^{*}\overline{%
\omega }^2$. Since $\phi $ has full rank the last equality implies that $%
dG=i_{\phi ^{*}IdH}\overline{\omega }^2$. From this follows that $I^{\prime
}dG=\phi _{*}IdH$, i.e. the gradient vector fields are $\phi $-related.

The system (\ref{dGdp}), (\ref{dGdq}) can be viewed as a system of linear
algebraic equations $Ah=g$ with the matrix

$$
A = \left( 
\begin{array}{cc}
\partial{x}/ \partial{p} & \partial 
{y} / \partial{p} \\ \partial{x}/ \partial{q} & \partial{y} / \partial{q} 
\end{array}
\right) 
$$

with the right hand part $g = [ \partial{\ \overline{G} } / \partial{p},
\partial{\overline{G}}/ \partial{q} ]^T$.

Both vectors $[(y,H)_{p,q},-(x,H)_{p,q}]^T$ and $[\partial{\overline{G}}/
\partial{x}, \partial{\overline{G}}/\partial{y}]^T$ satisfy the system: the
first one is the solution by assumption and the second one as the result of
the chain rule differentiation. The matrix $A$ is of a full rank and the
system is overdetermined. Hence the solution, if exists, is unique. This
proves that (\ref{dG}) holds under our assumptions, which is equivalent to
the $\phi$-relation of the corresponding vector fields. The proof is
complete.

Using formula $L_X = i_X \circ d + d \circ i_X$ for Lie derivative $L_X$
along vector field $X$, and the fact that $d\Omega^2=0$, we can get an
equivalent proposition:

\begin{proposition}
\label{proposition:HS_F_system1} A vector field $v$ on $M$ maps onto a
Hamiltonian field on $(N,\overline{\omega}^2)$ under surjective submersion $%
\phi: M \to N$ iff 2-form $\Omega^2 = \phi^* \overline{\omega}^2$ is an
invariant of 1-parametric Lie group generated by the vector field $v$.
\end{proposition}

Proof. Using infinitesimal criterion of the invariance we can conclude that $%
L_v\Omega ^2=0$. Next, expanding this using the formula $L_v=i_v\circ
d+d\circ i_v$, we get: $d\circ i_v\Omega ^2=0$ because $d\Omega ^2=0$ by
definition of symplectic form. This reduces Proposition \ref
{proposition:HS_F_system1} to Proposition \ref{proposition:HS_F_system}.

\newtheorem{corollary}{Corollary}

\begin{corollary}
\label{corollary:firstintegral} If $\phi : HS \to HS^{\prime}$ is morphism
of Hamiltonian systems then the function $\overline{G} = G \circ \phi$,
where $G$ is a Hamiltonian of $HS^{\prime}$, is the first integral of $HS$.
\end{corollary}

Proof. The corollary follows from the chain of equalities:

$$
(\overline{G},H)_{p,q} = IdH(d\overline{G}) = i_{IdH}\phi^*dG = (i_{IdH})^2
\phi^* \overline{\omega}^2 = 0 
$$

Here the first equality is by definition. The second one is an expansion of $%
d\overline{G}$. The third one is due to Proposition \ref
{proposition:HS_F_system1}. The last equality holds because of the skew
symmetry of 2-form $\phi ^{*}\overline{\omega }^2$, which is symplectic on $%
N $.

Remarks. In short this section says that we can reduce the order of a
Hamiltonian system by projecting a general vector field and then converting
the reduced dynamic system into the Hamiltonian form. The problem of
recognizing a Hamiltonian system in a general dynamic system was studied in
geometric mechanics (see \cite{koz} for linear quadratic case). The
Hamiltonian form always exists locally whenever we can present $2m-1$
independent first integrals, which is always possible in the vicinity of a
regular point \cite{birk}. We covered the subject in sufficient details
mainly to establish the framework for the following sections. Also note that
the propositions in this section can be easily generalized for the global
case, but we keep it local for consistency with the later discussion.

\section{Factorization of Optimal Control Systems}

Here we introduce a category of Lagrangian systems to set up a framework for
factorization of optimal control systems. Then we are going to establish a
connection between factorizations in categories of Lagrangian and
Hamiltonian systems.

\begin{definition}
\label{definition:LagrangianSystem} A Lagrangian system $LS$ is a triplet $%
(M, CDS, L)$, where $M$ is a manifold, $dimM=n$, $CDS$ is a controllable
dynamic system on $M$, and $L$ is a function $M \times U \to R$.
\end{definition}

By dynamical control system $CDS$ in the definition \ref
{definition:LagrangianSystem} we understand a system:

$$
\dot q^i = f^i(q,u), \qquad q \in M, i=\overline{1,n} 
$$

where the vector of control $u \in U=R^m$ is unconstrained.

A curve $\gamma : [0,T] \to R \times U$ (or, in coordinates, $%
\gamma(t)=(q(t),u(t)), t \in [0,T]$) is called a solution for $CDS$ if it
satisfies the equation $dq(t)/dt = f(q(t),u(t))$ for $\forall t \in [0,T]$.
Also we will call such curve admissible.

The function $L : M \times U \to R$ from the definition \ref
{definition:LagrangianSystem} defines a functional ${\cal L}(\gamma)$ on the
set of all admissible curves by the formula:

$$
{\cal L}(\gamma) = \int_{0}^{T} L(q(t),u(t))dt 
$$

{\it Solutions }of a Lagrangian system are the extremals of ${\cal L}(\gamma
)$ in the class of curves with fixed boundaries.

\begin{definition}
\label{definition:LS_solution} A solution $\gamma(t)=(q(t), u(t))$ of
Lagrangian system $LS=(M,CDS,L)$ is a curve providing a local extremum to
the functional ${\cal L}$ on the class of admissible curves with fixed
boundary points.
\end{definition}

Let $LS^{\prime}= (N, CDS^{\prime}, Q)$ be another Lagrangian system, such
that:\linebreak $dim N = \nu$, $CDS^{\prime}$ has the form $\dot y=F(y,v)$, $v \in V,
V = R^{\mu}$ and $Q=Q(y,v)$. Consider mapping $\Psi : M \times U \to N
\times V$, or in coordinates: $y=y(q,u), v=v(q,u)$.

\begin{definition}
\label{definition:LS_morphism} A mapping $\Psi : M \times U \to N \times V$
is called a morphism of Lagrangian systems from $LS$ to $LS^{\prime}$ if it
maps solutions of $LS$ into solutions of $LS^{\prime}$.
\end{definition}

In other words, if $\gamma =\gamma (t)$ is a solution of $LS$ then $\gamma
^{\prime }=\Psi \circ \gamma $, $\gamma ^{\prime
}(t)=(y(q(t),u(t)),v(q(t),u(t)))$ is a solution of $LS^{\prime }$. We will
denote a morphism of Lagrangian systems by the same mapping symbol $\Psi
:LS\to LS^{\prime }$. We will be interested in morphisms that are onto and
of full rank (factorizations) of Lagrangian systems.

Let $p_idq^i$ be the standard 1-form on $T^*M$.

\begin{definition}
\label{definition:PontryaginFunction} A function ${\cal H}(p,q,u) = p_i
f^i(q,u) - L(q,u)$ defined on $(T^*M) \times U$ is called Pontryagin
function of the Lagrangian system $LS = (M, CDS, L)$.
\end{definition}

We disregard singular systems, so $p_0$ in a more general Pontryagin
function $p_if^i(q,u)-p_0L(q,u)$ is never vanishing and we always have $%
p_0\equiv 1$.

\begin{definition}
\label{definition:NondegenerateLS} A Lagrangian system $LS$ is called
non-degenerate if it satisfies the following conditions:

\begin{enumerate}
\item  A system of nonlinear algebraic eqations $\partial {\cal H} / \partial%
{u^k} = 0, k=\overline{1,m}$ can be resolved with respect to $u$, the
soluition $\hat u = \hat h(p,q)$ is unique and the mapping $\hat u : T^*M
\to U$ is smooth.

\item  The matrix $f_u$ has full rank: 
\[
rank \biggl[ \frac{\partial{f^i}}{\partial{u^k}} \biggr] = m, 
\]
where $i=\overline{1,n}$ and $k=\overline{1,m}$.
\end{enumerate}
\end{definition}

The mapping $\hat u$ described above is the {\it optimal synthesis} for the
optimal control system.

Note that $m \leq n$ for a non-degenerate Lagrangian system (the dimension
of control space does not exceed the dimension of the phase space).

From here we will consider only non-degenerate Lagrangian systems.

The correspondence between control variables $u$ and dual variables $p$
established by optimal synthesis is not one to one in case $m<n$. Because of
that we need the following definition of observability.

\begin{definition}
\label{definition:ObservableInLSFunction} A function $S :T^*M \to R$ is
called observable in $LS$ if there exists a function $\Psi : M \times U \to R
$ such that the following diagram commutes

\begin{picture}(150,100)(-85,5) 
\put (20,80){$T^*M$}
\put (130,80){$M \times U$}
\put (85,15){$R$}
\put( 50, 85 ){\vector(1,0){70}}
\put( 35, 75 ){\vector(1,-1){45}}
\put( 145, 75 ){\vector(-1,-1){45}}
\put (85,90){$\Delta$}
\put (40,40){$S$}
\put (130,40){$\Psi$}
\end{picture} 

\begin{flushleft}
where $\Delta = \pi \times \hat{u}$ is morphism of fiber bundles
$\pi : T^*M \to M$ and $\pi' : M \times U \to M$. 
\end{flushleft}
\end{definition}

In other words, the definition requires that $S(p,q) = \Psi(q, \hat{u}(p,q))$
for observable in $LS$ function $S$. In the linear case this definition
corresponds to the observability defined in \cite{ar}.

The set of observable in $LS$ functions will be denoted as ${\cal F}_o(LS)$,
or ${\cal F}_o$ for brevity, when no confusion can happen.

We also need to define observability for morphisms of the Hamiltonian system
derived from a Lagrangian system. Let $HS^{\prime}=(T^*N, \overline{\omega}%
^2, G), dim N = \nu$ be a Hamiltonian system defined on $T^*N$ with a
natural symplectic form $\overline{\omega}^2=dx_i \wedge dy^i$ where $(x,y)$
are canonical coordinates on $T^*N$. Let $\phi:HS \to HS^{\prime}$ be
morphism. In coordinates $x=x(p,q), y= y(p,q)$. Let $\overline{\pi}:T^*N \to
N$ be a natural projection and $L_h$ be Lie derivative along vector field $%
h=IdH$ defined by $HS=(T^*M, \omega^2, H)$.

\begin{definition}
\label{definition:ObservableMorphism} Mapping $\phi: T^*M \to T^*N$ is
observable in $LS$ if functions $\overline{\pi} \circ \phi$ and $L_h%
\overline{\pi} \circ \phi$ are observable in $LS$.
\end{definition}

If $\phi :HS\to HS^{\prime }$ is a factorization of Hamiltonian systems and $%
\phi $ is observable in $LS$ then $HS^{\prime }$ is called an observable in $%
LS$ factorization of $HS$.

Also we will need a rather technical definition of a regular point that
would allow us to facilitate the proof of the main result later.

\begin{definition}
\label{definition:RegularPointOfLS} A point $(q_0,u_0) \in M \times U$ is
called a regular point of $LS$ if the set $R = \overline{\pi}^{-1}(q_0)
\times ker_{u_0}\Delta$ (in coordinates: $R=\{(p,q) \in T^*M: q=q_0,
u_0=\hat u (p,q_0) \}$) does not contain singular points of $HS$ and the
rank of the set of the functions $L_{h}^s q^i, i=\overline{0,n}, s = 0,1,
\dots$ is constant in the vicinity of each point in $R$.
\end{definition}

This type of regularity will turn out to be quite natural, but we will see
that only later in the discussion.

To show that the regular points do exist, consider Lagrangian system
corresponding to a linear-quadratic optimization problem. One can show that
regular points exist not only for linear-quadratic systems. Also it is
possible to somewhat relax requirements for the regular points but for the
price of much more technicalities that we would need to deal with. So we
presented a simpler but more restricting version of regularity.

The final preparation before formulating the main result of the theory of
factorization of Lagrangian systems is the following definition. A
Lagrangian system $LS^{\prime }$ is a {\it factor system} for $LS^{\prime }$
iff there exists morphism $\phi :LS\to LS^{\prime }$ which is a surjective
submersion.

\begin{proposition}
\label{theorem:FactorizationOfLS} Let $LS$ and $LS^{\prime}$ be Lagrangian
systems. Let $HS$ and $HS^{\prime}$ be corresponding Hamiltonian systems.
Then the following two conditions are equivalent:

\begin{enumerate}
\item  $LS^{\prime}$ is factorization of $LS$

\item  $HS^{\prime}$ is observable in $LS$ factorization of $HS$
\end{enumerate}
\end{proposition}

Proof. (2) $\Rightarrow$ (1). Given a morphism of Hamiltonian systems $\phi:
HS \to HS^{\prime}$ we will build the corresponding morphism $\psi : LS \to
LS^{\prime}$, which, by definition, is a mapping $\psi: M \times U \to N
\times V$ that maps extremals of $LS$ into extremals of $LS^{\prime}$.

The first half of the morphism mapping functions can be easily obtained from
the observability assumption: $y=y(p,q)=y(q,\hat{u}(p,q))$, since $y(p,q)$
is observable in $LS$. Hence on the extremals of $LS$ we have $y=y(q,u)$.

In coordinates, if $\gamma(t)=(p(t), q(t))$ is a solution of $HS$ then $%
\phi( \gamma(t))=(x(t),y(t))$ is a solution of $HS^{\prime}$ which
corresponds to an extremal $\tilde{\gamma}(t)=(y(t), \hat{v}(x(t),y(t))$ in $%
LS^{\prime}$ where $\hat{v}(x,y)$ is the optimal synthesis in $LS^{\prime}$.

From observability of $\phi$ we have $L_hy^i=\tilde{F}^i(q,\hat{u}(p,q))$
with some functions $\tilde{F}^i, i=\overline{1,\nu}$. Here $\hat{u}$ is
optimal synthesis in $LS$.

Since vector fields of $HS$ and $HS^{\prime}$ are $\phi$-related, we have
equalities: $\tilde{F^i}(q,\hat{u})=F^i(y(q,\hat{u}),\hat{v}), i=\overline{%
1,\nu}$, which hold on the trajectories of $HS$. The system $LS^{\prime}$ is
non-degenerate, thus we can resolve these equations. Indeed, these equations
are consistent on the trajectories of $HS$ with respect to $\hat{v}$:

\begin{eqnarray}
\hat{v}=v(q,u)  \label{synthesInFactorSystem}
\end{eqnarray}

It is easy to see that $y=y(p,u)$ and $v=v(q,u)$ are defining a morphism of
the Lagrangian systems. If $(q(t),u(t))$ is an extremal of $LS$, then there
exists $p(t)$ such that $(p(t),q(t))$ is a solution for $HS$. Morphism of
Hamiltonian systems maps this solution into a solution $(x(t),y(t))$ of $%
HS^{\prime}$. This solution defines an extremal $(y(t),v(t))$ with an
optimal synthesis $v(t)=\hat{v}(x(t),y(t))$. But on the trajectories of $HS$
holds (\ref{synthesInFactorSystem}), so $\hat{v}(x(t),y(t))=v(q(t), \hat{u}%
(p(t),q(t)))$ on the solutions of $HS$. This means that the mapping $v$
gives the same function of time as the optimal synthesis $\hat{v}$, so the
extremal of $LS$ was mapped into an extremal of $LS^{\prime}$.

(1) $\Rightarrow$ (2) Given morphism $\psi: LS \to LS^{\prime}$ of
Lagrangian systems we need to build morphism $\phi: T^*M \to T^*N$ of
Hamiltonian systems from $HS$ into $HS^{\prime}$ and show that $\phi$ is
observable in $LS$.

The first half of the morphism components is obvious: $y(p,q)= \linebreak y(q,\hat{u}%
(p,q))$ where $\hat{u}$ is the optimal synthesis and $y(q,u)$ is the first
part of the morphism $\psi: M \times U \to N \times V$ of the Lagrangian
systems. These functions are obviously observable.

If $\gamma(t)=(q(t),u(t))$ is an optimal trajectory in $LS$ then its image $%
\gamma^{\prime}=(y(q(t),u(t)), v(q(t),u(t)))$ is an extremal, hence an
admissible trajectory of $LS^{\prime}$. From this follows that on the
trajectories of $HS$ holds

$$
\frac{d}{dt}y(q,\hat{u})=F(y(q,\hat{u}), v(q,\hat{u})))) 
$$

Thus the functions $L_hy^i$ are observable. The observability of the
morphism is established and from now on we will write for brevity $%
y=y(p,q), \linebreak L_hy=\tilde F(p,q)$, collapsing the longer expression via $\hat
u(p,q)$.

Let $v=\hat v(x,y)$ be an optimal synthesis in $LS^{\prime }$. By assumption
the equality $F(y(p,q),\hat u(p,q))=\tilde F(p,q)$ holds whenever $%
(p(t),q(t))$ is a solution of $HS$. When this is the case, there exists a
function $x(t)$ such that together with $y(p(t),q(t))$ it satisfies $%
HS^{\prime }$. Under these conditions, our task is to find unknown
components $x(p,q)$ of the mapping $\phi :T^{*}M\to T^{*}N$ while we know
part of it $y(p,q)$ so that $\phi $ will be a morphism from $HS$ to $%
HS^{\prime }$. The Lemma from the next section claims that such functions $%
x(p,q)$ exist. Proving the Lemma will finish the proof of the Proposition.

\section{Existence of the Morphism}

The previous section left us with a partial mapping of Hamiltonian systems
that we need to extend to morphism. It turns out that the Hamiltonian
structure is not important for that so we will consider general dynamical
systems

\begin{eqnarray}
\dot x^j = \xi^j(x), j=\overline{1,m}  \label{Lemma:fieldX}
\end{eqnarray}

and

\begin{eqnarray}
\dot y^i = \eta^i(y,z), i=\overline{1,n_1}  \label{Lemma:fieldY1} \\
\dot z^i = \zeta^k(y,z), k=\overline{1,n_2}  \label{Lemma:fieldY2}
\end{eqnarray}

They define vector fields

\begin{eqnarray}
X = \xi ^j \frac{ \partial }{ \partial x ^ j} \\
Y = \eta ^i (y,z) \frac{ \partial }{ \partial y ^ i} + \zeta ^k (y,z) \frac{
\partial }{ \partial z ^ k}
\end{eqnarray}

\begin{flushleft}
where $j=\overline{1,m}, i=\overline{1,n_1}, k = \overline{1,n_2}$ 
and $ m \geq n = n_1+n_2 $
\end{flushleft}

The fields $X$ and $Y$ are defined in $U_1\in R^m$ and $U_2\in R^n$
respectively.

We assume that there exists a mapping $y=y(x)$ such that for each trajectory 
$x(t)$ of the field $X$ there exists a trajectory $(y(t),z(t))$ of the field 
$Y$ such that $y$ maps $x(t)$ onto the corresponding components of the
image. That means that if $x(t)$ is a solution of the system (\ref
{Lemma:fieldX}) and $y(t)=y(x(t))$ is provided by the mapping $y(x)$, then
there exists $z(t)$ such that $(y(x(t)),z(t))$ satisfies the system (\ref
{Lemma:fieldY1}, \ref{Lemma:fieldY2}).

As in Definition \ref{definition:RegularPointOfLS} we will call a point $x_0
\in U_1$ a regular one if $x_0$ is a regular point for the field $X$ and the
rank of the set of the functions $L^s_X \xi^i(x_0), s=0,1,..., i=\overline{%
1,m}$ is constant in the vicinity of $x_0$.

\newtheorem{lemma}{Lemma}

\begin{lemma}
\label{lemma:existenceofmorphism} Under the described assumptions in the
vicinity of the regular points of $X$ and $Y$ there exists a mapping $z(x)$,
not necessarily unique, such that the pair $(y(x),z(x))$ is a morphism of
dynamical systems from (\ref{Lemma:fieldX}) to (\ref{Lemma:fieldY1},\ref
{Lemma:fieldY2}).
\end{lemma}

The key observation for proving the Lemma is that the existence of a
complete mapping from (\ref{Lemma:fieldX}) to (\ref{Lemma:fieldY1}, \ref
{Lemma:fieldY2}) is equivalent to the existence of a solution of a PDE
system with {\it identical principal part}. This kind of systems was
explored by V.I.Elkin in \cite{el} and the proof of the Lemma relies on his
results.

Proof. Since for each solution $x(t)$ of (\ref{Lemma:fieldX}) there exists
some solution $(y(t),z(t))$ of (\ref{Lemma:fieldY1},\ref{Lemma:fieldY2}), we
can conclude that there exists a mapping $z=z(x)$ such that $F: x \to (y,z),
y=y(x), z= z(x)$ maps initial conditions $x_0$ of the solutions of (\ref
{Lemma:fieldX}) into initial conditions $(y_0, z_0)$ of the solutions of (%
\ref{Lemma:fieldY1},\ref{Lemma:fieldY2}).

At this point we can't claim yet that $F$ is the morphism we are looking for
since we need to show that it will differentiate properly along the field $X$
to map $X$ into $Y$. To show that we can differentiate $y(x)$ part of the
mapping along both fields. By the assumption the derivatives of $y$ have to
be the same along both fields since $y=y(x)$ maps solutions into a partial
solutions. That will result into a system of algebraic equations with
respect to $y^i, z^k$:

\begin{eqnarray}
y^i = y^i(x)  \label{Lemma1:algebraicSystem1} \\
L^s_Y y^i = L^s_X y^i(x)  \label{Lemma1:algebraicSystem2}
\end{eqnarray}

\begin{flushleft}
where $i=\overline{1,n_1}, k=\overline{1,n_2}$ and $s=1,2,...$.
\end{flushleft}

The graph of the mapping $F$ is a set $\{(x,y,z)\in U_1\times
U_2:y=y(x),z=z(x)\}$ contained inside of the manifold $M\in U_1\times U_2$,
which is defined by the equations (\ref{Lemma1:algebraicSystem1}, \ref
{Lemma1:algebraicSystem2}). The mapping $F$ is defined for all $x\in U_1$,
so for each $x$ there exists some solution $(y,z)$ of (\ref
{Lemma1:algebraicSystem1}, \ref{Lemma1:algebraicSystem2}). Thus dependent
equations in the system (\ref{Lemma1:algebraicSystem1}, \ref
{Lemma1:algebraicSystem2}) have to be identities with respect to $x$.

Next consider a regular point $(y_0,z_0)=(y(x_0),z(x_0))$ of the filed $Y$.
In the vicinity of the point $(x_0,y_0,z_0)\in U_1\times U_2$ the rank $r$
of the set of the functions

\begin{eqnarray}
L^s_Y y^i, i=\overline{1,n_1}, s=1,2,...  \label{Lemma:derivativesOfy}
\end{eqnarray}

is constant.

If the set is of the maximum rank $r=n_1+n_2$ then by the virtue of the
implicit function theorem the equations (\ref{Lemma1:algebraicSystem1}, \ref
{Lemma1:algebraicSystem2}) define a function $z=z(x)$ that together with $%
y=y(x)$ defines the required morphism of the dynamical systems.

To finish the proof we have to consider the case $r<n_1+n_2$. The fields $X$%
, $Y$ and $Z=X+Y$ have no singular points inside the area of consideration.
Thus they define a set of the first integrals. For the field $Z$ we will
have $n+m-1$ first integrals which are independent functions in the area.
The set of the integrals contain all the integrals $I^\nu(x), \nu=\overline{%
1,m-1}$ of the field $X$ and some functions $J^\alpha(x,y,z), \alpha=%
\overline{1,n}$, such that:

$$
det \biggl[ \frac{\partial J^\alpha}{dy^i} \: \bigg| 
\: \frac{\partial J^\alpha}{dz^j} \biggr ] \neq 0 
$$

where $i=\overline{1,n_1}, k=\overline{1,n_2}, n_1+n_2=n$.

We can add some of the functions $J^\alpha$ to the set (\ref
{Lemma:derivativesOfy}) to make it of full rank $n$. Without loosing any
generality we can assume that the functions used for that are $%
J^\beta(x,y,z), \beta=\overline{1,n-r}$. Let us add equations

\begin{eqnarray}
J^\beta(x,y,z) = I^\beta(x), \qquad \beta=\overline{1,n-r}
\label{Lemma1:AdditionalEqns}
\end{eqnarray}

to the system (\ref{Lemma1:algebraicSystem1},\ref{Lemma1:algebraicSystem2}).
The combined system (\ref{Lemma1:algebraicSystem1},\ref
{Lemma1:algebraicSystem2}, \ref{Lemma1:AdditionalEqns}) defines a manifold $%
M^{\prime }\in M$. If the manifold $M^{\prime \prime }$ is defined by the
equations (\ref{Lemma1:AdditionalEqns}), then $M^{\prime }=M\cap M^{\prime
\prime }$. It is easy to see that the field $Z$ is tangent to each of the
manifolds $M$ and $M^{\prime \prime }$, thus $Z$ is tangent to $M^{\prime }$%
. It is known from \cite{el} that such manifold defines smooth functions $%
y^i=y^i(x),z^k=z^k(x)$ satisfying a system of PDE with identical principal
part:

\begin{eqnarray}
\xi^j \frac{\partial y^i}{\partial x^j} = \eta^i(y,z) \\
\xi^j \frac{\partial z^k}{\partial x^j} = \zeta^k(y,z)
\end{eqnarray}

where $j=\overline{1,m}, i=\overline{1,n_1}, k=\overline{1,n_2}$. But this
system is exactly equivalent to the condition that the functions $%
y=y(x),z=z(x)$ define a morphism from (\ref{Lemma:fieldX}) to (\ref
{Lemma:fieldY1}, \ref{Lemma:fieldY2}). This concludes the proof.

Remark. Although the set (\ref{Lemma:derivativesOfy}) contains infinite
number of functions we can define the rank of this functions set with finite
number of differentiations and calculations of determinant. That allows to
reduce the system (\ref{Lemma1:algebraicSystem1}, \ref
{Lemma1:algebraicSystem2}) to a finite one that contains only independent
equations.

\section{Equations of Factorization. Building a Factorsystem.}

Let's denote ${\cal F}_1(LS)$ (or, in a short form, ${\cal F}_1$) the set of
the functions on $T^*M$ observable in $LS$ together with its first
derivative along $HS$. Obviously ${\cal F}_1 \subseteq {\cal F}_0$. Both
conditions $f \in {\cal F}_0$ and $f \in {\cal F}_1$ can be expressed in
terms of differential equations for $f$.

\begin{proposition}
\label{proposition:equationsOfFactorization} A Lagrangian system $LS$ has
non-trivial factorization iff there exist functions $y^i(p,q) \in {\cal F}_1
, \tilde{Q}(p,q) \in {\cal F}_0, x_i(p,q), i=\overline{1,\nu}, \nu < n$ such
that the functions $y^i(p,q), x_i(p,q), i=\overline{1,\nu}$ are independent
and the following equation is satisfied:

\begin{eqnarray}
L_{IdH}(x_i dy^i) = d\tilde{Q}  \label{EqnOfFactorization0}
\end{eqnarray}
\end{proposition}

To prepare the proof we will use Cartan formula for Lie derivative to
transform the equation (\ref{EqnOfFactorization0}) into its equivalent form:

\begin{eqnarray}
d (x_i(y^i,H)_{p,q} - \tilde{Q}) = i_{IdH}(dx_i \wedge dy^i)
\label{equationsOfFactorization}
\end{eqnarray}

Necessity. Let there exist a factorsystem $LS^{\prime}$

\begin{eqnarray}
\int_{0}^{T} Q(y,v) dt \rightarrow extr \\
\dot y = F(y,v)
\end{eqnarray}

where $y \in N, dim N = \nu$ and $v=\hat{v}(x,y)$ is its optimal synthesis.
Then the extremals of $LS^{\prime}$ are described by the equations:

\begin{eqnarray}
\dot y^i = \frac{\partial {\cal T}}{\partial x_i} \bigg| _{v=\hat{v}(x,y)} =
F^i(y,\hat{v}(x,y)) \\
\dot x_i = - \frac{\partial {\cal T}}{\partial{y^i}} \bigg| _{v=\hat{v}(x,y)}
\end{eqnarray}

Here ${\cal T}=x_iF^i(y,v)-\tilde Q(y,v)$ is Pontryagin function of the
Lagrangian system $LS^{\prime }$. By assumption we have the equality $%
F^i(y,\hat v(x,y))=(y^i,H)_{p,q}$. It follows from Proposition (\ref
{theorem:FactorizationOfLS}) that there exists a morphism $x=x(p,q),y=y(p,q)$
of the corresponding Hamiltonian systems. The Hamiltonian of the factor
system in coordinates $(p,q)$ is:

$$
\overline{G}(p,q) = x_i(p,q) (y^i(p,q),H)_{p,q} - Q(y(p,q),\hat{v}%
(x(p,q),y(p,q))) 
$$

From Proposition (\ref{proposition:HS_F_system}) we have equality $d%
\overline{G}=i_{IdH}(dx_i\wedge dy^i)$. It immediately leads to the equation
(\ref{equationsOfFactorization}), if we set: \\ $\tilde Q(p,q)= Q(y(p,q),\hat
v(x(p,q),y(p,q)))$.

It can be easily verified that $y^i \in {\cal F}_1$ and $\tilde{Q} \in {\cal %
F}_0$.

Sufficiency. We assume that functions $y^i(p,q) \in {\cal F}_1 , \tilde{Q}%
(p,q) \in {\cal F}_0, x_i(p,q),\linebreak i=\overline{1,\nu}, \nu < n$ satisfy
equation (\ref{EqnOfFactorization0}) and that functions $y^i(p,q), x_i(p,q),
i=\overline{1,\nu}$ are independent. Consider a function

$$
\overline{G} = x_i(p,q)(y^i(p,q),H)_{p,q} - \tilde{Q}(p,q) 
$$

and denote $(y^i(p,q),H)_{p,q} = \tilde{F}(p,q)$. We will show that $\tilde{F%
}=\overline{F} \circ \phi$ and $\tilde{Q}=\overline{Q} \circ \phi$, where $%
\phi$ is the mapping defined by $(x(p,q),y(p,q))$.

From the condition (\ref{equationsOfFactorization}) and Proposition (\ref
{proposition:HS_F_system}) follows that $\overline{G}=G \circ \phi$. Then,
according to Proposition (\ref{proposition:fieldsconnection}) on fields $%
\phi $-relation, we have:

\begin{eqnarray}
(y^i(p,q),H)_{p,q} = \frac{\partial G}{\partial x_i}(x(p,q), y(p,q)) = 
\overline{F}^i(x(p,q), y(p,q))  \label{Eq1.15}
\end{eqnarray}

Hence $\tilde F=\overline{F}\circ \phi $. Because of that we have $%
x_i(p,q)\tilde F^i(p,q)=\Theta \circ \phi $, where $\Theta =\Theta (x,y)=x_i%
\overline{F}^i(x,y)$. Finally, from $\tilde Q=$ $x_i(p,q)(y^i(p,q),H)_{p,q}-%
\overline{G}=\Theta \circ \phi -G\circ \phi $ we obtain that $\tilde Q=%
\overline{Q}\circ \phi $.

Next we will consider independent functions $v^k=\hat{v}^k(x,y), k=\overline{%
1,\mu}$ such that $\overline{F}^i = F^i(y,\hat{v}(x,y))$. The number $\mu$
of such functions can be derived from (\ref{Eq1.15}):

$$
rank \biggl[ \frac{\partial^2 G}{\partial x_i \partial x_j} \biggr] = rank 
\biggl[ \frac{\partial F^i}{\partial v^k} \frac{\hat{v}^k}{\partial{x_j}} 
\biggr]
$$

where $i,j=\overline{1,\nu}, k=\overline{1,\mu}$. Since

$$
rank \biggl[ \frac{\partial \hat{v}^k}{\partial x_j} \biggr] = \mu 
$$

then

$$
rank \biggl[ \frac{\partial F^i}{\partial v^k} \biggr] = \mu 
$$

Having $\overline{F}^i=F^i(y,\hat{v})$ we can show that $\overline{Q}$
depends on $x$ via $\hat{v}$: $\overline{Q}=Q(y,\hat{v})$. From the
definition of $\overline{Q}$ we see that $\overline{Q}=x_i\overline{F}^i(y,%
\hat{v}) - \overline{G}(x,y)$. Differentiation of $\overline{Q}$ reveals
linear dependence of the gradients:

$$
\frac{\partial \overline{Q}}{\partial x_j} = F^j(y,\hat{v}) + x_i \frac{%
\partial F^i}{\partial v^k} \frac{\partial \hat{v}^k}{\partial x_j} - \frac{%
\partial \overline{G}}{\partial x_j} = (x_i \frac{\partial F^i}{\partial v^k}%
) \frac{\partial \hat{v}^k}{\partial x_j} 
$$

Here the last equality follows from (\ref{Eq1.15}). From the linear
dependence we conclude that $\overline{Q}(x,y)=\overline{Q}(y,\hat{v}(x,y))$.

Now we can form a Lagrangian system:

\begin{eqnarray}
\int^T_0 Q(y,v) dt \rightarrow extr \\
\dot y^i = F^i(y,v), \qquad i=\overline{1,\nu}
\end{eqnarray}

We will show that $v=\hat{v}(x,y)$ is the optimal synthesis for this system.
The Pontryagin function for this system is: ${\cal T}(x,y,v) = x_iF^i(y,v) -
Q(y,v)$. By construction ${\cal T}(x,y,\hat{v}(y,v)) = G(x,y)$, thus

$$
\frac{\partial G}{\partial x_j} = \frac{\partial {\cal T}}{\partial x_j}%
(x,y,v) \bigg|_{v=\hat{v}(x,y)} + \frac{\partial {\cal T}}{\partial v^k}%
(x,y,v) \bigg|_{v=\hat{v}(x,y)} \frac{\partial \hat{v}^k}{\partial x_j} 
$$

On the other hand

$$
\frac{\partial {\cal T}}{\partial x_j}(x,y,v)\bigg|_{v=\hat{v}(x,y)} = \frac{%
\partial G}{\partial x_j} = F^j(y,\hat{v}(x,y)) 
$$

Comparing last two equalities we see that:

$$
\frac{\partial {\cal T}}{\partial v^k}(x,y,v) \bigg|_{v=\hat{v}(x,y)} \frac{%
\partial \hat{v}^k}{\partial x_j} = 0 
$$

Since $\hat{v}^k$ are independent functions of $x$, we can satisfy the last
equality only if

$$
\frac{\partial {\cal T}}{\partial v^k}(x,y,v) \bigg|_{v=\hat{v}(x,y)} = 0 
$$

But this is exactly the definition of the optimal synthesis as a stationary
point of Pontryagin function.

From the observability of $\tilde{F}(p,q)$ follows observability of the
synthesis $\hat{v}^k(x(p,q),y(p,q))$. To proof this consider:

$$
\frac{\partial \tilde{F}^j}{\partial p_s} = A^j_i(p,q) \frac{\partial \hat{u}%
^a}{\partial p_s} 
$$

with some functions $A^j_i(p,q)$. Their existence follows from the
observability of $\tilde{F}(p,q)$. On the other hand

$$
\frac{\partial \tilde{F}^j}{\partial p_s} = \frac{\partial \tilde{F}^j}{%
\partial y^i} \frac{\partial y^i}{\partial p_s} + \frac{\partial \tilde{F}^j%
}{\partial \hat{v}^k} \frac{\partial \hat{v}^k}{\partial p_s} = \frac{%
\partial \tilde{F}^j}{\partial y^i} B^i_a(p,q) \frac{\partial \hat{u}^a}{%
\partial p_s} + \frac{\partial \tilde{F}^j}{\partial \hat{v}^k} \frac{%
\partial \hat{v}^k}{\partial p_s} 
$$

where again the existence of the functions $B_a^i(p,q)$ follows from the
chain rule for differentiation of a compound function and from observability
of $y^j(p,q)$ and their derivatives. Combining all together we see that the
gradient $\partial \hat v^k/\partial p$ is linearly expressed via the
gradients $\partial \hat u^a/\partial p$ since the matrix $[\partial
F^j/\partial u^a]$ is of the maximum rank. This proves observability of $%
\hat v$ in $LS$.

Thus we built a factor system from the solution of (\ref{EqnOfFactorization0}%
) and also we built a morphism into the factor system. This concludes the
proof.

We will call the equations from the Proposition (\ref
{proposition:equationsOfFactorization}) (and their equivalents (\ref
{equationsOfFactorization})) "equations of factorization".

Note that calculating the outer derivative $d$ on both sides of (\ref
{equationsOfFactorization}) we get $L_{IdH}(dx^i\wedge dy_i)=0$ which is in
agreement with Proposition \ref{proposition:HS_F_system}. Also it is easy to
verify that the identity morphism provides a solution for the equations of
factorization(\ref{equationsOfFactorization}): $x=p,y=q,\tilde Q=L$, which
transforms these equations into the identity:

$$
d(p_i\frac{\partial H}{\partial p_i}-L)=dH 
$$

Remark. A quite standard note is that systems with explicit dependence on
time $t$ can be reduced to the investigated stationary case. We can add a
new equation $\dot t = 1$ and a new pair of boundary conditions for the new
variable $t(0)=0$ and $t(T)=T$. The only difference of the new variable from
the rest is that its boundary conditions are always the same. This does not
affect our reasoning which was for the fixed ends case anyway. Finally, a
field of extremals remains a local object in this case too.

{\bf Example 1.} An optimal control system

\begin{eqnarray}
\int_0^T (q_1u_1u_2 + q_1q_2) dt \rightarrow extr  \label{Example1a} \\
\dot q_1 = u_1, \qquad \dot q_2 = u_2  \label{Example1b}
\end{eqnarray}

has a factorization. This was shown by A.N.Chernoplekov in \cite{ch}. The
system (\ref{Example1a}, \ref{Example1b}) is especially well suited to show
that the factorization theory developed here is quite natural and
generalizes factorization of variational problems in \cite{ch}.

The Hamiltonian equations for the optimal control system (\ref{Example1a},%
\ref{Example1b}) are:

\begin{eqnarray}
\dot p_1 = \frac{p_1p_2}{q_1^2} + q_2, \qquad \dot q_1= \frac{p_2}{q_1}
\label{Example1.1HamsystemA} \\
\dot p_2 = q_1, \qquad \dot q_2 = \frac{p_1}{q_1}
\label{Example1.1HamsystemB}
\end{eqnarray}

The mapping $(x=2p_2, y=q_1^2)$ defines a morphism into the Hamiltonian
factor system

$$
\dot x = 2 \sqrt{y}, \qquad \dot y = x 
$$

The Lagrangian system that corresponds to this Hamiltonian system is

\begin{eqnarray}
\int_0^T(\frac{1}{2}v^2+\frac{4}{3}y^{3/2})dt \rightarrow extr \\
\dot y = v
\end{eqnarray}

and the optimal synthesis is $\hat{v}=x$. The conditions of observability in
this case are satisfied automatically. The morphism of the Lagrangian
systems in this case is

$$
y=q_1^2, \qquad v=2q_1u_1 
$$

Also let's write down the equations of the factorization and check that they
are satisfied:

$$
(2p_2d(q_1^2))^. = d( \frac{1}{2}(2p_2)^2 + \frac{4}{3}q_1^3) 
$$

Here $()^{.}$ denotes differentiation $L_{IdH}$ along the vector field of
the original Hamiltonian system (\ref{Example1.1HamsystemA},\ref
{Example1.1HamsystemB}).

The morphism of the Hamiltonian systems in this example is essentially
non-symplectic, meaning that there is no coordinate change that will make it
symplectic. It is easily follows from the observation that for any two
functions $f_1(q_1,p_2), f_2(q_1,p_2)$ we have $(f_1,f_2)_{p,q}=0$.

It is also possible to verify that the Lagrangian system from this example
does not allow any symmetries as defined in \cite{ech}. To show that we will
search for a vector filed in the form (see proof of Theorem 1 in \cite{ech}):

$$
\tilde{\xi} = \xi_1(q_1, q_2) \frac{\partial}{\partial q_1} + \xi_2(q_1,
q_2) \frac{\partial}{\partial q_2} + \zeta_1(q_1,q_2,u_1,u_2) \frac{\partial%
}{\partial u_1} + \zeta_2(q_1,q_2,u_1,u_2) \frac{\partial}{\partial u_2} 
$$

If $X$ is vector field defined by dynamical system (\ref{Example1b}), then
invariance of Lagrangian system with respect to a field $\tilde{\xi}$ is
given in \cite{ech} by conditions $L_{\tilde{\xi}}X=0$ (invariance of the
vector field $X$) and $L_{\tilde{\xi}}L=0$ (invariance of the Lagrangian $%
L(q,u)$).

Expanding these conditions into the system of PDE we obtain:

\begin{eqnarray}
\zeta_1 = u_1 \frac{\partial \xi_1}{\partial q_1} + u_2 \frac{\partial \xi_1}{\partial q_2} \\
\zeta_2 = u_1 \frac{\partial \xi_2}{\partial q_1} + u_2 \frac{\partial \xi_2}{\partial q_2} \\
\xi_1(u_1u_2+q_2) + \xi_2 q_1+q_1 u_2 \zeta_1+ q_1 u_1 \zeta_2 = 0
\end{eqnarray}

We can substitute $\zeta_i$ into the last equation. Then we can break it
into a system of equations by powers of $u_{1,2}$ since the solution $\xi$
does not depend on $u$. The resulting system

\begin{eqnarray}
\xi_1q_2 + \xi_2q_1 = 0 \\
q_1 \frac{\partial\xi_2}{\partial q_1} = 0, \qquad q_1 \frac{\partial\xi_1}{\partial q_2} = 0 \\
\xi_1 + q_1 \frac{\partial\xi_1}{\partial q_1} + q_1 \frac{\partial\xi_2}{\partial q_2} = 0
\end{eqnarray}

has only trivial solution.

Thus the Lagrangian system in this example has no symmetries in terms of 
\cite{ech} yet it allows order reduction within the introduced category of
Lagrangian systems.

{\bf Example 2}. A system

\begin{eqnarray}
\int_0^T(q_1q_2 + \frac{1}{2}u_1^2)dt \rightarrow extr \\
\dot q_1 = q_2-u_1, \qquad \dot q_2 = q_1 + u_1
\end{eqnarray}

offers another example of order reduction via factorization in the category
of Lagrangian systems.

One of the possible strategies of order reduction for this system is to sum
up equations of the control system and eliminate $u_1$. Then we can
integrate the resulting ODE, but that will result in introducing time $t$ in
the right part.

However solving the equations of factorization will allow for more elegant
order reduction. The solution that leads to a simpler system is

$$
x=q_1-q_2+p_2-p_1, \qquad y = p_2-p_1, \qquad \tilde{Q} = \frac{1}{2}%
(q_1-q_2+p_2-p_1)^2 - (p_2-p_1)^2 
$$

The corresponding factor-system is:

\begin{eqnarray}
\int_0^T \frac{1}{2}(v^2-y^2) dt \rightarrow extr \\
\dot y = v
\end{eqnarray}

with the morphism

$$
y=u_1 \qquad v = q_1-q_2+u_1 
$$

into the factorsystem.

{\bf Example 3}. The so called ''horizontal decomposition'' (see \cite{pa1})
can be achieved for the system:

\begin{eqnarray}
\int_0^T(\frac{1}{2}u_2^2 - u_1u_3 - \frac{1}{4}q_1^2 - \frac{1}{4}q_3^2) dt
\rightarrow extr \\
\dot q_1 = u_1, \qquad \dot q_2 = q_1 + u_2 + u_3, \qquad \dot q_3 = q_2 -
u_2 + u_3
\end{eqnarray}

This system is equivalent to the pair of two independent factor systems:

\begin{eqnarray}
\int_0^T(y_1y_2 + \frac{1}{2}v_1^2) dt \rightarrow extr \\
\dot y_1 = y_2, \qquad \dot y_2 = v_1
\end{eqnarray}

and

\begin{eqnarray}
\int_0^T\frac{1}{2}(v_3^2 - y_3^2) dt \rightarrow extr \\
\dot y_3=v_3
\end{eqnarray}

The corresponding morphisms of the Lagrangian systems are:

\begin{eqnarray}
y_1 = q_1, \qquad y_2 = u_1, \qquad v_1 = \frac{1}{2}(q_2+q_3) \\
y_3 = -u_2 \qquad v_3 = \frac{1}{2}(q_2-q_3)
\end{eqnarray}

And the morphisms of their Hamiltonian systems are:

\begin{eqnarray}
x_1 = p_1, \qquad x_2 = \frac{1}{2}(q_2+q_3), \qquad y_1=q_1, \qquad y_2 =
-p_2 - p_3 \\
x_3 = \frac{1}{2}(q_2-q_3), \qquad y_3 = p_3 - p_2
\end{eqnarray}

\section{Boundary Conditions,\\ Constrained Control and Factorization}

In this section we will discuss how boundary conditions transform under
factorization of Lagrangian systems.

We will say that the Lagrangian system defines an optimal control problem if
we specify and fix some boundary conditions allowing to pose a fixed ends
boundary conditions problem for the corresponding Hamiltonian system. As an
example we will consider conditions $q(0)=q_0,q(T)=q_1$. Under the
factorization mapping they transform into boundary conditions in factor
spaces $T^{*}N$ and $N\times V$. The transformed boundary conditions define
some manifolds where the trajectory of the factor system has to start and to
end. Depending on the dimension of these manifolds we can end up with over-,
well- or under- determined boundary problem for factorsystem.

The over- and well-determinied cases are directly useful. If the boundary
problem for the original system has a solution it is obvious that the factor
problem is also solvable, even if it may appear overdetermined. Such
factorization allows for classical hierarchial control when we can solve the
factor problem and then extend its solution to the solution of the original
problem. The discussion of such well-behaving factorization naturally falls
into the framework developed in {Y.N.Pavlovskii, V.I.Elkin}~\cite{pa1} and
will be essentially the same as in the paper~\cite{ch} by {A.N.Chernoplekov}
for the case of variational systems.

The case of under-determined boundary factor problem is less obvious.

{\bf Example 4}. Consider Lagrangian system

\begin{eqnarray}
\int_0^T(u_1u_2 + q_2) dt \rightarrow extr \\
\dot q_1 = u_1, \qquad \dot q_2 = u_2
\end{eqnarray}

It has factor system

\begin{eqnarray}
\int_0^T \frac{1}{2}v^2 dt \rightarrow extr \\
\dot y = v
\end{eqnarray}

with the morphism $y=u_1u_2,v=u_2$ into it. The fixed boundary conditions $%
q_1(0),q_2(0),q_1(T),q_2(T)$ define a unique trajectory in the original
Lagrangian system. But, when mapped into the factorsystem, they don't
provide enough information to build the appropriate boundary conditions for
the factor system. In fact they don't impose any restrictions on the
trajectory ends at all since the boundary manifolds coincide with the entire
space.

The case of an under-determined factor system does not allow to immediately
benefit from knowing the factorsystem. We still can build a hierarchy out of
factor- and quotient- systems. Its functioning may rely on a differential
game with two players: the center (factorsystem) and subordinate (quotient)
system. The goals of functioning of both players are the same but the
natural information structure does not allow to achieve the optimum without
using some additional interaction between the players. A differential games
resulting from such factorization may resemble situations when the center
operates in terms that are not directly related to the reality in which the
subordinate system has to function, even though they share the same goal.
Thus the practical value of such factorization is less obvious.

Note that the discussed here factorization of Lagrangian systems can be
easily generalized to an optimal control problem with free end(s). That
follows from the locality of the field of extremals. Instead of the boundary
conditions $q(0)=q_0,q(T)=q_1$ we will have $q(0)=q_0,p(T)=0$. That change
of boundary conditions does not affect any of our constructions.

Finally we will briefly touch the case of constrained controls. Since the
approach we used here is based on smooth objects, any direct generalization
to the constrained case might me problematic. However in many practical
cases we can approximate the original constrained optimal control problem
with a smooth unconstrained one by introducing smooth penalty functions to
represent constraints. Apparently it could be done in many different ways
potentially leading to different factorizations or no factorizations at all.

\section{Conclusion}

In this paper we interpreted order reduction of optimal control systems as
factorization in the category of Lagrangian systems. We established
sufficient and necessary conditions of factorization for Lagrangian systems.
Factorization can be described in terms of the corresponding Hamiltonian
systems that appear from Pontryagin's maximum principle.

Morphisms of Hamiltonian systems that we use in the paper differ from the
classical. Our definition does not require the mapping into factor system to
be symplectic. That means that we are not necessarily able to extend it to a
canonical change of coordinates. Because of that Hamiltonian factor systems
does not preserve the original symplectic form on the base space. Instead
the symplectic form in the factor space has to be invariant under the flow
of the original Hamiltonian system. This invariance of the symplectic form
gives sufficient and necessary condition of the factorization of Hamiltonian
systems. Naturally, Hamiltonian of the factor system turns out to be first
integral of the original Hamiltonian system. Finally, factorization of
Hamiltonian system allows to build Lagrangian factor system iff the morphism
is an observable mapping. Observability here means that the mapping and its
Lie derivative along the original Hamiltonian field depends on the dual
variables only via optimal synthesis functions.

We discussed factorization of boundary conditions for Lagrangian systems.
The mapping of the fixed boundary conditions under the factorization does
not always allow to obtain a well defined boundary condition problem for the
factor systems. Also we discussed some of the possible interpretations of
that situation from the point of view of differential games.

\section{Acknowledgment}

I'd like to thank my scientific supervisor Y.N.Pavlovskii for his guidance
and for my initial interest in the problem. Also I would like to thank
A.N.Chernoplekov, V.I.Elkin, A.P.Krischenko, S.A.Kutepov, T.G.Smirnova, and 
G.N.Yakovenko for their feedback and valuable discussions.

\end{document}